\documentclass[12pt, twoside, eqno]{article}
\usepackage{latexsym}
\usepackage{amssymb}
\usepackage{amsfonts}
\begin{document}

\noindent {\bf \large On hypersemigroups}\bigskip

\medskip

\noindent{\bf Niovi Kehayopulu}\bigskip

\noindent November 27, 2015

\bigskip

\noindent{\small
{\bf Abstract.} We prove that a nonempty subset $B$ of a regular 
hypersemigroup $H$ is a bi-ideal of $H$ if and only if it is 
represented in the form $B=A*C$ where $A$ is a right ideal and $C$ a 
left ideal of $H$. We also show that an hypersemigroup $H$ is regular 
if and only if the right and the left ideals of $H$ are idempotent, 
and for every right ideal $A$ and every left ideal $B$ of $H$, the 
product $A*B$ is a quasi-ideal of $H$. Our aim is not just to add a 
publication on hypersemigroups but, mainly, to publish a paper which 
serves as an example to show what an hypersemigroup is and give the 
right information concerning this structure. We never work directly 
on an hypersemigroup. If we want to get a result on an 
hypersemigroup, then we have to prove it first for a semigroup and 
transfer its proof to hypersemigroup. But there is further 
interesting information concerning this structure as well, we will 
deal with at another time.\medskip

\medskip

 }

\section {Introduction}
Many results on semigroups based on ideals, bi-ideals and 
quasi-ideals are due to S. Lajos. In his paper in [2] he proved that 
if $S$ is a regular semigroup, then a nonempty subset $B$ of $S$ is a 
bi-ideal of $S$ if and only if there exists a right ideal $R$ and a 
left ideal $L$ of $S$ such that $B=RL$. J. Calais proved in [1] that 
a semigroup $S$ is regular if and only if the right and the left 
ideals of $S$ are idempotent and for every right ideal $A$ and every 
left ideal $B$ of $S$, the product $AB$ is a quasi-ideal of $S$. In  
the present note we present analogous results for hypersemigroups. We 
tried to use sets instead of elements in an attempt to show that, 
exactly as in semigroups, for the results on hypersemigroups based on 
ideals, points do not play any essential role, but the sets which 
shows their pointless character.

An {\it hypergroupoid} is a nonempty set $H$ with an hyperoperation 
$$\circ : H\times H \rightarrow {\cal P}^*(H) \mid (a,b) \rightarrow 
a\circ b$$on $H$ and an operation $$* : {\cal P}^*(H)\times {\cal 
P}^*(H) \rightarrow {\cal P}^*(H) \mid (A,B) \rightarrow A*B$$ on 
${\cal P}^*(H)$ (induced by the operation of $H$) such that 
$$A*B=\bigcup\limits_{(a,b) \in\,A\times B} {(a\circ b)}$$ for every 
$A,B\in {\cal P}^*(H)$ (${\cal P}^*(H)$ denotes the set of nonempty 
subsets of $H$).\smallskip

An hypergroupoid $(H,\circ,*)$ is called {\it hypersemigroup} if 
$$\{x\}*(y\circ z)=(x\circ y)*\{z\}$$for all $x,y,z\in H$. Since 
$\{x\}*\{y\}=x\circ y$ for every $x,y\in H$, the associativity of $H$ 
can be also defined by 
$$\{x\}*{\Big(}\{y\}*\{z\}{\Big)}={\Big(}\{x\}*\{y\}{\Big)}*\{z\}$$ 
for every $x,y,z\in H$. \smallskip

A nonempty subset $A$ of an hypergroupoid $H$ is called a {\it left} 
(resp. {\it right}) ideal of $H$ if $H*A\subseteq A$ (resp. 
$A*H\subseteq A)$. A subset of $H$ which is both a left and a right 
ideal of $H$ is called an {\it ideal} of $H$. For any nonempty subset 
$A$ of an hypersemigroup $H$, we denote by $R(A)$, $L(A)$ and $I(A)$ 
the right ideal, left ideal and the ideal of $H$, respectively, 
generated by $A$ and we have

$R(A)=A\cup (A*H)$, $L(A)=A\cup (H*A)$, and

$I(A)=A\cup (H*A)\cup (A*H)\cup (H*A*H)$.\\A nonempty subset $A$ of 
an hypergroupoid $H$ is called {\it idempotent} if $A=A*A$.
\section{Main results}
By the definition of the hypergroupoid we have the following lemma 
which, though clear, plays an essential role in the investigation: 
\medskip

\noindent{\bf Lemma 1.} {\it If $H$ is an hypergroupoid and $A,B\in 
{\cal P}^*(H)$, then

$(1)$ $x\in A*B$ $\Longleftrightarrow$ $x\in a\circ b$ for some $a\in 
A$, $b\in B$.

$(2)$ If $a\in A$ and $b\in B$, then $a\circ b\subseteq 
A*B$.}\medskip

\noindent{\bf Lemma 2.} {\it If $H$ is an hypergroupoid then, for any 
$A,B,C,D\in {\cal P}^*(H)$, we have

$(1)$ $A\subseteq B
\mbox { and } C\subseteq D\;\Longrightarrow\; A*C\subseteq B*D$.

$(2)$ $A\subseteq B \Longrightarrow A*C\subseteq B*C \mbox { and } 
C*A\subseteq C*B.$

$(3)$ $H*H\subseteq H$.

$(4)$ $H*A\subseteq H$ and $A*H\subseteq H$.}\medskip

\noindent {\bf Proof.} (1) Let $A\subseteq B$ and $C\subseteq D$. If 
$x\in A*C$, then $x\in a\circ b$ for some $a\in A$, $b\in B$, then 
$x\in a\circ b$ for some $a\in B$, $b\in D$, so $x\in B*D$. The proof 
of the rest is similar. $\hfill\Box$\medskip

\noindent{\bf Lemma 3.} {\it If $H$ be an hypergroupoid and $A_i,B\in 
{\cal P}^*(H)$, $i\in I$, then
$$(\bigcup\limits_{i \in I} {{A_i}} )*B = \bigcup\limits_{i \in I} 
{({A_i}} *B).$$}{\bf Proof.} Let $x\in (\bigcup\limits_{i \in I} 
{{A_i}} )*B$. Then $x\in a\circ b$ for some $a\in (\bigcup\limits_{i 
\in I} {{A_i}} )$, $b\in B$. Since $a\in A_j$ for some $j\in I$ and 
$b\in B$, we have $a\circ b\subseteq A_j*B\subseteq \bigcup\limits_{i 
\in I} {({A_i}} *B)$.
Let now $x\in A_j*B$ for some $j\in I$. Then $x\in a\circ b$ for some 
$a\in A_j$, $b\in B$. Since $a\in A_j\subseteq \bigcup\limits_{i \in 
I} {{A_i}} $ and $b\in B$, we have $a\circ b\subseteq 
(\bigcup\limits_{i \in I} {{A_i}} )*B$. Then we get $x\in 
(\bigcup\limits_{i \in I} {{A_i}} )*B$. $\hfill\Box$\medskip

\noindent{\bf Lemma 4.} {\it If $H$ is an hypersemigroup, then the 
operation ``$*$" on ${\cal P^*}(H)$ is associative, that is
$(A*B)*C=A*(B*C)$ for any any $A,B,C\in {\cal P}^*(H)$.}\medskip

\noindent According to this lemma, instead of $(A*B)*C$ or $A*(B*C)$, 
we write $A*B*C$. The proof of this lemma is due to M. Tsingelis and 
it is based on Lemma 3. \medskip

\noindent{\bf Definition 5.} An hypersemigroup $H$ is called {\it 
regular} if, for any $A\in {\cal P}^*(H)$, we have
$A\subseteq A*H*A$.\medskip

\noindent{\bf Definition 6.} Let $H$ be an hypersemigroup. A nonempty 
subset $B$ of $H$ is called a  {\it bi-ideal} of $H$ if 
$$B*H*B\subseteq B.$${\bf Proposition 7.} {\it Let H be an 
hypersemigroup. If C is a right ideal of $H$ and $D\in {\cal P}^*(H)$ 
(or D a left ideal of H and $C\in {\cal P}^*(H)$), then the set 
$B=C*D$ is a bi-ideal of H}.\medskip

\noindent{\bf Proof.} Let $C$ be a right ideal of $H$ and $D$ a 
nonempty subset of $H$ such that $B=C*D$. Then $B*H*B\subseteq B$. 
Indeed: We 
have\begin{eqnarray*}B*H*B&=&(C*D)*H*(C*D)=C*(D*H*C)*D\\&\subseteq& 
(C*H)*D\subseteq C*D=B.\end{eqnarray*} $\hfill\Box$

\noindent{\bf Theorem 8.} {\it If H is a regular hypersemigroup then, 
for every bi-ideal B of H there exists a right ideal C and a left 
ideal  D of H such that $$B=C*D.$$}{\bf Proof.} Let $B$ be a bi-ideal 
of $H$. Then $B*H*B\subseteq B.$ Since $H$ is regular, we have 
$B\subseteq B*H*B$, thus we have $B=B*H*B$. On the other hand, we 
have\begin{eqnarray*}R(B)*L(B)&=&{\Big(}B\cup (B*H){\Big)}* 
{\Big(}B\cup (H*B){\Big)}\\&=&(B*B)\cup {\Big(}(B*H)*B{\Big)}\cup 
{\Big(}B*(H*B){\Big)}\cup {\Big(}(B*H)*(H*B){\Big)}\\&=&(B*B)\cup 
(B*H*B)\cup {\Big(}B*(H*H)*B{\Big)}.\end{eqnarray*}Since 
$H*H\subseteq H$, we have $B*(H*H)\subseteq B*H$ and 
${\Big(}B*(H*H){\Big)}*B\subseteq (B*H)*B$. Thus we 
have$$R(B)*L(B)=(B*B)\cup (B*H*B)=(B*B)\cup B.$$In addition, 
$B*B=(B*H*B)*B={\Big(}B*(H*B){\Big)}*B$. Since $H*B\subseteq H$, we 
get $B*(H*B)\subseteq B*H$ and ${\Big(}B*(H*B){\Big)}*B\subseteq 
(B*H)*B$. Thus we obtain $B*B\subseteq B*H*B=B$, and $R(B)*L(B)=B$, 
where $R(B)$ is a right ideal and $L(B)$ is a left ideal of $H$. 
$\hfill\Box$\\By Proposition 7 and Theorem 8, we have the following 
theorem\medskip

\noindent{\bf Theorem 9.} {\it If H is a regular hypersemigroup then 
B is a bi-ideal of H if and only if there exists a right ideal C and 
a left ideal D of H such that $B=C*D$.} \medskip

\noindent{\bf Definition 10.} If $H$ is an hypergroupoid, a nonempty 
subset $Q$ of $H$ is called a {\it quasi-ideal} of $H$ if $$(Q*H)\cap 
(H*Q)\subseteq Q.$${\bf Lemma 11.} {\it If H is an hypergroupoid, A a 
right ideal and B a left ideal of H, then $A\cap B\in{\cal 
P}^*(H)$.}\medskip

\noindent{\bf Proof.} Take an element $a\in A$ and an element $b\in 
B$ $(A,B\not=\emptyset)$. Since $\{a\}\subseteq A$ and 
$\{b\}\subseteq B$, we have $\{a\}*\{b\}\subseteq A*B$. Since
$\{a\}*\{b\}=a\circ b$, we have $a\circ b\subseteq A*B\subseteq 
A*H\subseteq A$ and $a\circ b\subseteq A*B\subseteq H*B\subseteq B$, 
so $a\circ b\subseteq A\cap B$. Since $a\circ b\in{\cal P}^*(H)$, we 
have $A\cap B\in {\cal P}^*(H)$. $\hfill\Box$\medskip

\noindent{\bf Theorem 12.} {\it An hypersemigroup H is regular if and 
only if the right and the left ideals of H are idempotent, and for 
every right ideal A and every left ideal B of H, the product $A*B$ is 
a quasi-ideal of $H$.}\medskip

\noindent{\bf Proof.} $\Longrightarrow$. Let $A$ be a right ideal of 
$H$. Since $H$ is regular, we have $$A\subseteq (A*H)*A\subseteq 
A*A\subseteq A*H\subseteq A,$$ thus we have $A*A=A$. Similarly, the 
left ideals of $H$ are idempotent. Let now $A$ be a right ideal and 
$B$ a left ideal of $H$. Then $A\cap B$ is a quasi-ideal of $H$. In 
fact: Since $A\cap B\subseteq A,B$, we have $(A\cap B)*H\subseteq 
A*H\subseteq A$ and $H*(A\cap B)\subseteq H*B\subseteq B$, thus we 
have $${\Big (}(A\cap B)*H{\Big)}\cap {\Big (}H*(A\cap 
B){\Big)}\subseteq A\cap B \;\;\;\;\;\;(*)$$ On the other hand, 
$A\cap B=A*B$. Indeed: Since $A\cap B\in{\cal P}^*(H)$ and $H$ is 
regular, we have
$$A\cap B\subseteq (A\cap B)*H*(A\cap B)\subseteq 
A*H*B=(A*H)*B\subseteq A*B.$$ We also have $A*B\subseteq A*H\subseteq 
A$ and $A*B\subseteq H*B\subseteq B$, so $A*B\subseteq A\cap B$, then 
$A*B=A\cap B$. Thus the set $A*B$ is a quasi-ideal of $H$.\smallskip

\noindent $\Longleftarrow$. Let $A$ be a nonempty subset of $H$. By 
hypothesis, we have\begin{eqnarray*}A&\subseteq& 
R(A)=R(A)*R(A)={\Big(}A\cup (A*H){\Big)}*{\Big(}A\cup 
(A*H){\Big)}\\&=&(A*A)\cup {\Big(}(A*H)*A{\Big)}\cup 
{\Big(}A*(A*H){\Big)}\cup {\Big(}(A*H)*(A*H){\Big)}\\&=&(A*A)\cup 
{\Big(}A*(H*A){\Big)}\cup {\Big(}A*(A*H){\Big)}\cup 
{\Big(}A*(H*A)*H{\Big)}.
\end{eqnarray*}Since $A\subseteq H$, we get $A*A\subseteq A*H$. Since 
$H*A\subseteq H$, we have $A*(H*A)\subseteq A*H$. Since $A*H\subseteq 
H$,
$A*(A*H)\subseteq A*H$. Since $A*(H*A)\subseteq A*H$, we have 
$A*(H*A)*H\subseteq A*(H*H)\subseteq A*H$. Thus we have $A\subseteq 
A*H$. Since $L(A)$ is a left ideal of $H$, in a similar way we get 
$A\subseteq H*A$. Therefore we obtain $$A\subseteq (A*H)\cap (H*A).$$ 
Since $A*H$ is a right ideal and $H*A$ is a left ideal of $H$, by 
hypothesis, they are idempotent, and we 
have\begin{eqnarray*}(A*H)\cap (H*A)&=&{\Big(}(A*H)*(A*H){\Big)}\cap 
{\Big(}(H*A)* (H*A){\Big)}\\&=&{\Big(}(A*H*A)*H{\Big)}\cap 
{\Big(}H*(A*H*A){\Big)}.\end{eqnarray*}On the other hand, $A*H*A$ is 
a quasi-ideal of $H$. Indeed: Since $H$ is a right ideal of $H$, by 
hypothesis it is idempotent that is, $H*H=H$. Thus we 
have$$A*H*A=A*(H*H)*A=(A*H)*(H*A).$$Since $A*H$ is a right ideal and 
$H*A$ a left ideal of $H$, by hypothesis, $(A*H)*(H*A)$ is a 
quasi-ideal of $H$. So $A*H*A$ is a quasi-ideal of $H$ as well, which 
means that$${\Big(}(A*H*A)*H{\Big)}\cap 
{\Big(}H*(A*H*A){\Big)}\subseteq A*H*A.$$Then we have $A\subseteq 
A*H*A$, and $H$ is regular. $\hfill\Box$\medskip

\noindent{\bf Corollary 13.} [2] {\it If S is a regular semigroup, 
then B is a bi-ideal of S if and only if there exists a right ideal C 
and a left ideal D of S such that $B=CD.$}\medskip

\noindent{\bf Corollary 14.} [1] {\it A semigroup S is regular if and 
only if the right and the left ideals of S are idempotent, and for 
every right ideal A and every left ideal B of S, the product $AB$ is 
a quasi-ideal of $S$.}\bigskip

\noindent{\bf Note.} Although the corresponding results on semigroups 
can be also obtained as application of the results of this paper, 
exactly as in $\Gamma$-semigroups, we never work directly on 
hypersemigroups. We first have to prove the results on semigroups and 
then to transfer them to hypersemigroups. So although the results on 
hypersemigroups generalize the corresponding results on semigroups, 
they are completely based on the results on semigroups. We can do 
nothing on an hypersemigroup if we do not have it first for a 
semigroup. We wrote this paper in an attempt to show the way we pass 
from semigroups to hypersemigroups and give the right information 
concerning this structure. We will give further information 
concerning this structure in a next paper.{\small
\bigskip

\noindent University of Athens, Department of Mathematics,
15784 Panepistimiopolis, Greece\\
email: nkehayop@math.uoa.gr

\end{document}